\documentclass[11pt]{article}
\usepackage{amsmath}
\usepackage{amssymb}
\usepackage{amscd}
\usepackage{latexsym}  

\numberwithin{equation}{section}

\newtheorem{propo}[equation]{Proposition}

\newtheorem{theor}[equation]{Theorem}
\newtheorem{lem}[equation]{Lemma}

\newtheorem{coro}[equation]{Corollary}
\newcommand{\qed}{\hfill $\Box$\\}
\newenvironment{pf}{{\it Proof.}}{\qed}
\newtheorem{acknowledgment}{Acknowledgment}

\begin{document}
\title{\bf Groups that do not act by automorphisms of codimension-one
foliations}
\author{R. Feres and D. Witte}
 \footnotetext{Department of Mathematics (CB 1146), Washington University,
 St. Louis, MO 63130, USA.}
 \footnotetext{Department of Mathematics, Oklahoma State University,
 Stillwater, OK 74078, USA.}

\date{July 29, 2000}
\maketitle
\begin{abstract}
 Let $\Gamma$ be a finitely generated group having the property that any
action of any finite-index
subgroup of~$\Gamma$ by homeomorphisms of the circle must have a finite
orbit. (By a theorem of
\'E.~Ghys,   lattices in   simple Lie groups of real rank at least $2$ have
this property.) Suppose
that   such a $\Gamma$ acts on a compact manifold $M$
 by automorphisms of  a codimension-one $C^2$ foliation, $\cal F$. We show
that if $\cal F$ has a
compact leaf, then some finite-index subgroup of~$\Gamma$ fixes a compact
leaf of~$\cal F$.
Furthermore, we give sufficient conditions for some finite-index subgroup
of $\Gamma$
to fix each leaf of $\cal F$.
 \end{abstract}

\section{Introduction and statement of results}
$M$ will denote a compact, connected,  boundaryless, smooth manifold of
dimension $n$. Let $\cal F$
be a $C^r$ foliation of $M$ by smooth leaves, $r\geq 2$. It will be assumed
that $\cal F$ is a
transversely oriented,   codimension-one foliation.

Let  ${\cal D}^s(M,{\cal F}) $  denote the group of $C^s$ automorphisms of
$\cal F$, $s\geq 0$; that is, the group
of $C^s$ diffeomorphisms of $M$ that map leaves to leaves.
The normal subgroup consisting of automorphisms of $\cal F$ that send each
leaf to itself will be
denoted ${\cal D}^s(M,{\cal F})_0$, and will be called the group of {\em
inner automorphisms} of
the foliation. The quotient $${\cal O}^s(M,{\cal F})={\cal D}^s(M,{\cal
F})/{\cal D}^s(M,{\cal F})_0 $$
will be called the group of {\em transverse automorphisms} (or {\em outer
automorphisms}) of $\cal F$.
When a group $\Gamma$ acts by automorphisms of $\cal F$ so as to define a
homomorphism into ${\cal D}^s(M,{\cal F})$,
  the action will be called {\em transversely finite} if $\Gamma$ projects to
a finite subgroup of ${\cal O}^s(M,{\cal F})$.

 The general question that motivates the results of the present  paper can
be stated thus: Given
$(M,{\cal F})$, what groups can act by  automorphisms of $\cal F$ so that
the action is not transversely finite,
and what groups cannot?  For example, if $\cal F$
is the foliation by fibers of a product manifold $M=B\times L$, with leaves
$\{b\}\times L$, $b\in B$,
then any group that acts nonfinitely on $B$ by $C^s$-diffeomorphisms also
acts nonfinitely by
outer-automorphisms of $\cal F$ (e.g.,
by setting the action on $L$ to be trivial). In this case, ${\cal
O}^s(M,{\cal F})$ is the group
of $C^s$-diffeomorphisms of $B$.  If, on the other hand, $M=\Bbb T^n$ (the
$n$-dimensional flat torus)
and $\cal F$ is the foliation  by planes parallel to an irrational
hyperplane $F\subset \Bbb R^n$,
then it is an elementary fact that
 ${\cal O}^s(M,{\cal F})$ is isomorphic to $H:= (\Gamma\ltimes {\Bbb
R}^n)/(\Gamma_0\ltimes F)$, where
$\Gamma$ is the stabilizer of $F$ in $GL(n,\Bbb  Z)$ and $\Gamma_0$ is the
subgroup of $\Gamma$ that
acts trivially on the quotient $\Bbb R^n/F$. In this case, only groups that
admit homomorphisms with
nonfinite image in $H$ can have nonfinite actions by outer-automorphisms of
$\cal F$.
(Allowing big codimension, it is quite easy to construct, say,
topologically transitive foliated bundles, with large groups of smooth
outer-automorphisms.)

A more specific question that will be addressed here is the following:
Suppose that no (topological, say) action of
a group $\Gamma$ on the circle
 yields nontrivial dynamics (that is, non-finite action). Does $\Gamma$
admit nontrivial ``transverse dynamics'' on some codimension-one
foliation of a compact manifold? The two theorems given below provide
support for the negative answer.

\begin{theor}\label{atomthm}
 Suppose that  $\Gamma$ is a finitely generated discrete group such that every
homomorphism from a finite-index subgroup of $\Gamma$
  into the group of homeomorphisms of the circle has a finite orbit on the
circle.
Also suppose  that $\Gamma$ acts by $C^0$-automorphisms of  $(M,{\cal F})$.

If $\cal F$ has a closed leaf, then some closed leaf of $\cal F$ is fixed
by some
subgroup of finite index in $\Gamma$.
\end{theor}

\begin{coro}\label{noatomcor}
 Suppose that  $\Gamma$ is a finitely generated discrete group such that every
homomorphism from a finite-index subgroup of~$\Gamma$
  into the group of homeomorphisms of the circle has a finite orbit on the
circle.
Also suppose  that $\Gamma$ acts by $C^0$-automorphisms of  $(M,{\cal F})$.

If $\cal F$ admits a bounded transverse invariant measure, then it also
admits a transverse
invariant measure $\mu$ that is    invariant under $\Gamma$,  and  every
leaf  in the support of
$\mu$ is sent to itself under the action of  a  finite-index
subgroup~$\Gamma'$ of~$\Gamma$.
 \end{coro}

The next theorem provides a class of foliations on which the
$\Gamma$-action is transversely finite.
A foliation is said to be  {\em almost without holonomy} if the germinal
holonomy
groups of all the non-compact leaves are trivial \cite[IV-2.11,
p.~251]{Godbillon}.

\begin{theor}\label{bigtwo}
Suppose that  $\Gamma$ is a  finitely generated discrete group such that
every homomorphism
from a finite-index subgroup of $\Gamma$ into the group of homeomorphisms
of the circle has a finite
orbit on the circle.
 Also suppose   that
 $\cal F$ is almost without  holonomy. Then every homomorphism of $\Gamma$
into ${\cal D}^1(M,{\cal F})$
yields a transversely finite action. In particular, the conclusion holds if
the non-compact leaves of $\cal F$
are simply connected.
\end{theor}

The following theorem of \'E.~Ghys \cite{ghys} provides examples of groups
that satisfy the
 requirements  of Theorems \ref{atomthm} and \ref{bigtwo}, that is, groups
for which
every homomorphism  into the group of homeomorphisms of the circle has a
finite orbit.
Most of these examples were also established by M.~Burger and N.~Monod
\cite{burger, burger2}.

\begin{theor}[{{\rm Ghys \cite{ghys}}}]\label{ghys}
Suppose that  $\Gamma$ is an irreducible lattice in a connected,
semisimple, real Lie group $G$ of
 real rank at least $2$, and  that there is no continuous homomorphism from
$G$ onto
$PSL(2,\Bbb R)$. Then every homomorphism from $\Gamma$ into the group of
homeomorphisms of the circle, $\Bbb T^1$, has a finite orbit. Furthermore,
if $\Gamma$ acts   by
$C^1$ diffeomorphisms, then some finite-index subgroup of~$\Gamma$ acts
trivially on~$\Bbb T^1$.
 \end{theor}

\begin{acknowledgment}
 {\upshape
 The authors would like to thank the Isaac Newton Institute for
Mathematical Sciences (Cambridge, U.K.) for providing a congenial environment
 to carry out this research; and
  S. Hurder and L. Conlon, for a number of helpful discussions and guidance
in the literature of
codimension-one foliations.
 D.W.\ was partially supported by a grant from the National Science
Foundation (DMS-9801136).
 }
 \end{acknowledgment}

\section{General facts about codimension-one foliations}

 We use \cite{conlon} as our source for basic facts about codimension-one
foliations. These results
were first proved by R.~Sacksteder, P.~Dippolito, G.~Hector, A.~Haefliger,
and others. Some of the
relevant papers are \cite{dipo,hec,mmt,sacks}.

The foliation $\cal F$ will be said to be {\em without holonomy} if the
germinal holonomy group
of each leaf of $\cal F$ is trivial.

Let $\cal L$ denote a smooth one-dimensional foliation of $M$ everywhere
transverse to $\cal F$.
(Cf.\ \cite[5.1.2]{conlon}.) It is convenient to work with {\em biregular}
coordinate charts for
the pair $({\cal F}, {\cal L})$. These are charts  that  define foliation
boxes for
both $\cal L$ and $\cal F$ simultaneously, having   local coordinate maps
  $\varphi\colon U\subset M\rightarrow V\subset \Bbb R^{n-1}\times \Bbb R^1$,
 $\varphi(p)= \bigl( x(p),y(p) \bigr)$, such that $x=\text{constant}$
corresponds to plaques of $\cal L$
while $y=\text{constant}$ corresponds to plaques of $\cal F$.
A {\em biregular cover} is an atlas comprised of biregular coordinate
charts. Such covers exist.
(Cf.\ \cite[5.1.4]{conlon}.) From now on $\cal L$ will denote   a fixed
transverse foliation to $\cal F$ and
any foliation box will be assumed without mention to be biregular.

An open $\cal F$-saturated set $U$ is a called a {\em foliated product} if
it is connected
and ${\cal L}|_U$ fibers $U$ by open intervals over some
$(n-1)$-dimensional manifold $B$.
Since $\cal F$ is orientable, a foliated product  is a trivial interval
bundle, homeomorphic to
 $B\times (0,1)$ (although the foliation need not be the product
foliation). Each leaf of $\cal F$ in $U$
with the restriction to it of the bundle map
is a covering space of $B$. We note, in particular, that each closed
transversal that meets
$U$ has to meet every leaf in $U$.
Let $d$ be the topological metric on $U$ induced by the restriction to $U$
of a Riemannian metric on $M$
and denote by $\widehat{U}$ the completion of $U$ in the metric $d$.

An $\cal F$-saturated set $U$ is a called a {\em foliated bundle} if it is
connected
and ${\cal L}|_U$ fibers $U$ over some $(n-1)$-dimensional manifold $B$.
(This is more general than
a foliated product, because there is no restriction on the fibers.) An
$\cal F$-saturated set $U$ is a called
a {\em trivially foliated product} if there is a connected $1$-manifold~$F$
(possibly with boundary), a connected
$(n-1)$-dimensional manifold $B$, and a diffeomorphism from~$U$ to $B
\times F$, that carries ${\cal L}|_U$ and
${\cal F}|_U$ to the product foliations of $B \times F$.

\begin{theor}[{{\rm Dippolito \cite{dipo}}}]\label{lemgen}
Let $U$ be a connected $\cal F$-saturated open set.
$\widehat{U}$ is a connected manifold with finitely many boundary
components. The interior of
$\widehat{U}$ is $U$ and the inclusion $i\colon U\rightarrow M$ extends  to
an immersion $\hat{\iota}$ of $\widehat{U}$
into $M$ that sends the boundary components of $\widehat{U}$ onto boundary
leaves of $U$. If $L'$ is a boundary leaf
of $U$ then $\hat{\iota}^{-1}(L')$ consists of one or two  components of
the boundary of $\widehat{U}$, each
component being mapped bijectively to $L'$ by $\hat{\iota}$. Both $\cal F$
and $\cal L$ pull-back under
$\hat{\iota}$ to well-defined foliations on $\widehat{U}$. If $U$ is  a
foliated product, then $\widehat{U}$
is a foliated bundle whose fibers are compact intervals.
\end{theor}

 \begin{pf}
 This is \cite[5.2.10, 5.2.11, 5.2.12]{conlon} as well as the remarks after
5.2.12 of the same
reference.
 \end{pf}

The foliation of $\widehat{U}$ obtained by the pull-back of $\cal F$ will
be denoted $\widehat{\cal F}$.

\begin{theor}[{{\rm Sacksteder \cite{sacks}}}]\label{general-one}
 Let $\cal F$ be a transversely orientable foliation of class $C^2$ and
codimension one on a compact
manifold. Then the following are equivalent:
 \begin{enumerate}
 \item There exists a bounded transverse invariant measure  $\mu$.
 \item Either $\cal  F$ has a compact leaf or $\cal F$ is without holonomy.
\end{enumerate}
\end{theor}

 \begin{pf}
 This is \cite[2.3.8]{hector}.
 \end{pf}

\begin{theor}\label{general-two}
 Let $\cal F$ be a transversely orientable foliation of class $C^2$ and
codimension-one of a compact
manifold $M$. Let $U$ be a connected $\cal F$-saturated open set and
 suppose that ${\cal F}|_U$ is without holonomy. Then either every leaf of
${\cal F}|_U$ is closed
and ${\cal F}|_U$ fibers over a connected $1$-manifold, or each leaf of
${\cal F}|_U$ is dense in
$U$. Furthermore, if ${\cal F}|_U$ is a fibration over a $1$-manifold $B$, but
$(\widehat{U},\widehat{\cal F})$ is not a trivially foliated product, then
$B\cong S^1$.
 \end{theor}

\begin{pf}
 This is \cite[9.1.4, 9.1.6]{conlon}.
 \end{pf}

\begin{theor}\label{general-three}
 Let $\cal F$ be a transversely orientable foliation of class $C^2$ and
codimension-one. Let
$U\subset M$ be  a connected, nonempty, open, $\cal F$-saturated set.
 \begin{enumerate}
 \item \label{partone} Suppose
that $(U, {\cal F}|_{U})$ is without holonomy. Then  there is a $C^0$ flow
$$\Phi\colon \Bbb R\times \widehat{U}\rightarrow \widehat{U}$$ that fixes
the points of $\partial\widehat{U}$,
 carries leaves diffeomorphically to leaves and is transitive on the set of
leaves of ${\cal F}|_U$.
Furthermore, $\cal F|_U$ admits a transverse-invariant  nonatomic  measure
$\mu$ of
full support that is bounded on compact subsets of $U$ and    assigns to
each transverse arc
$\{\Phi_t(p) | t\in (a,b)\}$ the measure $b-a$.

 \item Conversely, if $(U,{\cal F}|_U)$ admits a nonatomic transverse
invariant measure $\mu$
of full support, then $\cal F|_U$ is without holonomy and there exists a
continuous flow $\Phi$ that carries
leaves diffeomorphically to leaves, transitively on the set of leaves.
$\Phi$ is related to $\mu$ in the
way    described in part~\ref{partone} of this  theorem.
 \end{enumerate}
 \end{theor}

\begin{pf}
 This is a special case of \cite[9.2.1]{conlon}. The transverse-invariant
measure is described in the
proof given in the reference.  See also \cite{mmt}.
 \end{pf}

\section{Proof of Theorem~\ref{atomthm}} \label{withatoms}

Suppose that a $\Gamma$-action  satisfying the assumptions of
Theorem~\ref{atomthm} has been fixed.
We may assume, by passing to a finite-index subgroup of~$\Gamma$, that the
$\Gamma$-action preserves
the transverse orientation for~$\cal F$.

The $\Gamma$-action will be said to be {\em fixing} if there exists  a
finite-index
subgroup~$\Gamma'$ of~$\Gamma$ and a compact leaf~$L$ such that every
element of
$\Gamma'$  sends $L$ to itself.

A nonempty $\cal F$-saturated subset $\cal P$ of $M$ will be called $\cal
F$-{\em perfect} if
for every differentiable curve $\alpha\colon (-a,a)\rightarrow M$, $a>0$,
transverse to $\cal F$, the
intersection of $\cal P$ with the image of $\alpha$ is a perfect set.

\begin{lem}\label{lemmaP}
 If $\cal F$ has a compact leaf and the $\Gamma$-action is not
fixing, then there exists a compact, $\cal F$-saturated, $\cal F$-perfect,
$\Gamma$-invariant
 set ${\cal P}\subset M$,
 that is the union
of compact, mutually homeomorphic  leaves.
 \end{lem}

\begin{pf}
  Let $\cal C$ denote the union of all compact leaves of $\cal F$ of a same
homeomorphism type.
A nonempty set of this kind exists since $\cal F$ has a compact leaf. It is
also clear that $\cal C$
is invariant under every automorphism of~$\cal F$. Furthermore, by a
theorem of Haefliger,
\cite[6.1.1]{conlon},
 $\cal C$ is compact.
Let $${\cal A}=\{A\subset {\cal C} \mid A \text{ is compact, nonempty,
$\cal F$-saturated, and
$\Gamma$-invariant}\}.$$ By Zorn's lemma, $\cal A$ has an element $A$ that
is minimal under
inclusion.

We define the {\em derived  set} $A'$ of $A$ as the subset of $A$ comprised
of the union of all leaves
$L\subset A$ whose points are limits of sequences in the complement of $L$
in $A$.
If $A$ is the union of finitely many leaves, a finite-index subgroup
of~$\Gamma$ would
send each of those finitely many  leaves to itself,
contradicting the assumption that the $\Gamma$-action is not fixing.
Therefore $A'$ is nonempty.
It is easy to see that $A'\in \cal A$, so $A=A'$ by the minimality of $A$.
Therefore, ${\cal P}:=A$ satisfies the properties required in the lemma.
\end{pf}

\begin{propo}\label{PisM}
 If ${\cal P}=M$, then the foliation fibers over the circle.
\end{propo}

\begin{pf}
This is immediate from Theorem \ref{general-one} and the easy fact that
$\cal F$ is, in this case, without holonomy.
\end{pf}

 The connected components of   $M - {\cal P}$ will be called the {\em gaps}
of $\cal P$.

\begin{lem}\label{notfix}
 Suppose that ${\cal P}$ is a proper subset  of $M$, and
the $\Gamma$-action is not fixing. Then
 \begin{enumerate}
 \item each gap is bounded by two leaves of~$\cal P$, and
 \item there exists a finite open cover $\{V_1,\dots, V_k\}$ of $M$ by
foliated products such that
each $V_i$ is bounded by two   leaves in $\cal P$.
 \end{enumerate}
\end{lem}

\begin{pf}
We first remark that
each connected component of $W=M-{\cal P}$   is a foliated product.
In fact, by  \cite[5.2.9]{conlon},  only finitely many connected components
of $W$ fail  to be  foliated products.
Suppose that there are connected components of $W$ which are not foliated
products and denote them by
  $W_1,\dots, W_l$. Since any homeomorphism
   $\gamma\in \Gamma$
must send each $W_i$ into  some (possibly the same) $W_j$, one obtains a
homomorphism
of $\Gamma$ into the group of permutations of $l$ symbols, from which it
follows that some
finite-index subgroup of~$\Gamma$ sends each $W_i$ to itself. In particular,
this subgroup would permute the boundary components of a $W_i$. By
\cite[5.2.5]{conlon} the boundary
of a connected $\cal F$-saturated open set consists of the union of a
finite number of leaves,
which in this case must be elements of $\cal P$. But then a finite-index
subgroup of~$\Gamma$ would
send one leaf in $\cal P$ to itself, contradicting the assumption that the
action is not fixing.

A leaf in $\cal P$ will be called a {\em border leaf} of $\cal P$ if it is
a component
of the boundary of a gap. If $\cal P$ is not all
of $M$, there must be a countable infinity of gaps, since otherwise $\cal
P$ would be contained
in the union of  the finitely many boundary leaves  of a finite number of
connected $\cal F$-saturated open sets. Each border leaf $L$  of $\cal P$
is a boundary component of a gap
and on the side of $L$ opposite the gap
a sequence of   leaves in $\cal P$ accumulates on $L$.

If a leaf $L$ of $\cal P$ does not bound a gap, then $L$ is a limit of
sequences of leaves in $\cal P$ on
both of its sides, so that
  \cite[5.3.4]{conlon} (due to Dippolito \cite{dipo})   immediately yields
a foliated product neighborhood of $L$.

We claim that boundary of each gap of $\cal P$ consists of two (distinct)
leaves in $\cal P$, and
that the closure of each gap is contained in    a foliated product bounded
by two leaves in $\cal P$.
 The interiors of these foliated products together with the interiors of
the foliated products of
the previous paragraph form an open cover  for $M$. Since $M$ is compact,
we can extract a finite
open cover.

 All that is left is to prove the claim.
Let $W$ denote a gap of $\cal P$. We have seen that  it
is a foliated product. Denote by $B$ the base manifold.
 The boundary of $W$ consists of two (distinct) leaves in $\cal P$. (There
are not more than
two leaves by Theorem~\ref{lemgen}. If a single
leaf of $\cal P$ bounded $W$ on both sides, this would be an isolated leaf,
which is not the case.)
The boundary leaves are homeomorphic to $B$ (notice that
the boundary leaves of $\widehat{W}$ are homeomorphic to the
base of the fibration on $\widehat{W}$) and $\widehat{W}$ maps bijectively
onto the closure of $W$ under
the map $\hat{\iota}$.
Since the boundary leaves $L_1$ and $L_2$ are compact, it is possible to
find   neighborhoods $U_1$ and $U_2$
of $L_1$ and $L_2$, respectively, such that $U_i|_{\cal L}$ is a trivial
bundle over $L_i$, $i=1,2$,
not necessarily $\cal F$-saturated. To obtain    $\cal F$-saturated $U_i$,
one applies
\cite[5.3.4]{conlon}. The union of $W$, $U_1$ and $U_2$ (for sufficiently
small $U_i$) gives the desired  neighborhood.
\end{pf}

For $x \in {\cal P}$, let $L_x$ be the leaf of~$\cal F$ that contains~$x$.
Define an equivalence
relation~$\sim$ on~$\cal P$ by specifying that $x \sim y$ if either
 $L_x \cup L_y$ is the boundary of a gap of~$\cal P$, or  $L_x = L_y$. (In
particular, if ${\cal P} = M$, then $x \sim y$ if and only if $L_x = L_y$.)
Note that each
equivalence class is either a leaf or the union of two leaves.

\begin{lem}\label{col}
 If $\Gamma$ is as in Theorem~\ref{atomthm} and $\cal F$ has a compact
leaf, but the
$\Gamma$-action is not fixing, then ${\cal P}/{\sim}$ is homeomorphic to~$S^1$.
 \end{lem}

\begin{pf}
 We may assume ${\cal P}$ is a proper subset of~$M$; otherwise, the desired
conclusion follows
from Lemma~\ref{PisM}.
 It is immediate from Lemma~\ref{notfix} (and the ``waterfall
construction'' described in
\cite[3.3.7]{conlon}) that  there exists a closed transversal, $\alpha$, of
$\cal F$ that intersects
each leaf of $\cal P$ exactly once.  The intersection   $\alpha\cap \cal P$
is a perfect set in the
embedded circle $\alpha$, and the saturations of the  gaps of this perfect
set are the gaps of $\cal
P$. Thus, the desired conclusion follows from the elementary observation
that, by identifying the two
endpoints of each of the  gaps of $\alpha\cap \cal P$ to a single point, we
obtain a quotient that
is homeomorphic to a circle.
 \end{pf}

Suppose the $\Gamma$-action is not fixing, and let $\cal P$ be as in
Lemma~\ref{lemmaP}. The action
of~$\Gamma$ on~$\cal P$ factors through to an action of~$\Gamma$ by
homeomorphisms of~${\cal
P}/{\sim}$. Now Lemma~\ref{col} implies that ${\cal P}/{\sim}$ is
homeomorphic to~$S^1$, so, by
assumption, $\Gamma$ must have a finite orbit on~${\cal P}/{\sim}$. This
finite orbit yields a
$\Gamma$-invariant, finite collection of compact leaves in~${\cal P}$. Then
some finite-index
subgroup of~$\Gamma$ fixes each of these compact leaves. This proves
Theorem~\ref{atomthm}.

\section{Proof of Corollary~\ref{noatomcor}} \label{noatoms}

Suppose that $\cal F$ admits a transverse invariant measure.

If $\cal F$ has a compact leaf, then Theorem~\ref{atomthm} implies that
some finite-index
subgroup~$\Gamma'$ of~$\Gamma$ fixes some compact leaf~$L$. Then $\Gamma'$
fixes the atomic
measure~$\mu$ supported on the single leaf~$L$, so the conclusion of
Corollary~\ref{noatomcor} holds.

Thus, we may assume that no leaf of $\cal F$ is compact. Therefore, from
Theorems~\ref{general-one}
and~\ref{general-two}, it follows that $\cal F$ is a minimal foliation
(that is, every leaf is dense
in $M$). Then Corollary~\ref{corobigtwo} below completes the proof of
Corollary~\ref{noatomcor}.

Corollary~\ref{corobigtwo} is stated in greater generality than needed for
the proof of
Corollary~\ref{noatomcor} because it  will be used in the given form for
the proof of Theorem
\ref{bigtwo}.

\begin{propo}[{{\rm cf.\ \cite[Thm.~X.2.3.3, p.~272]{hector}}}] \label{mu'=cmu}
 Let $U$ be a connected, $\Gamma$-invariant, $\cal F$-saturated open set.
We suppose that the
boundary of $U$ is either empty or consists of finitely many
$\Gamma$-invariant compact leaves. We
also suppose that each leaf of $\cal F$ in $U$ is dense in $U$. Let $\mu$
be a transverse invariant
measure on $U$ that is bounded on each compact subset of~$U$.

If $\mu'$ is another transverse invariant
measure on $U$ that is bounded on each compact subset of~$U$, then $\mu'$
is a scalar multiple
of~$\mu$.
 \end{propo}

\begin{coro} \label{corobigtwo}
 Let $U$ be a connected, $\Gamma$-invariant, $\cal F$-saturated open set.
We suppose that the
boundary of $U$ is either empty or consists of finitely many
$\Gamma$-invariant compact leaves. We
also suppose that each leaf of $\cal F$ in $U$ is dense in $U$. Let $\mu$
be a transverse invariant
measure on $U$ which is bounded on any compact subset of $U$. Then
 \begin{enumerate}
 \item \label{corobigtwo-mu} $\mu$ is $\Gamma$-invariant; and
 \item \label{corobigtwo-leaf} $[\Gamma,\Gamma]$ is a finite-index subgroup
of~$\Gamma$ that fixes
each leaf in~$U$.
 \end{enumerate}
 \end{coro}

\begin{pf} (\ref{corobigtwo-mu}) For each $\gamma \in \Gamma$,
Proposition~\ref{mu'=cmu} implies
there is some $c(\gamma) \in \Bbb R^+$, such that $\gamma_* \mu = c(\gamma)
\, \mu$. It is easy to
see that $c \colon \Gamma \to \Bbb R^+$ is a homomorphism, so, because
$\Bbb R^+$ is abelian and has
no nontrivial finite subgroups, Lemma~\ref{finitecommutator} implies
$c(\Gamma) = 1$. Thus, $\mu$ is
$\Gamma$-invariant.

(\ref{corobigtwo-leaf}) We use the notations of Theorem
\ref{general-three}, where the Sacksteder
flow $\Phi_t$ is defined.
 Integration of $ \mu$
over closed curves representing elements   of $\pi_1(M)$ yields a homomorphism
$\rho\colon \pi_1(M)\rightarrow \Bbb R$ whose image group, $P(\mu)$, is
called the {\em group of periods}
of $\mu$ \cite[9.3.4]{conlon}. This is a finitely generated abelian
subgroup of $\Bbb R$. The group of periods can
be characterized by the following property (\cite[9.3.6]{conlon}): $\Phi_t$
sends every leaf to
itself exactly when $t\in P(\mu)$.

Define for each $p\in M$ and each $\gamma\in \Gamma$ a class $[t]\in \Bbb
R/P(\mu)$ where $t$
is any real number such that $\Phi_t(p)$ lies in the leaf  of $\gamma(p)$.
Then the correspondence
$\gamma\mapsto [t]$ gives a well-defined homomorphism, $h$ from $\Gamma$
into $\Bbb R/P(\mu)$.
Because $\Bbb R/P(\mu)$ is abelian, we know that $[\Gamma,\Gamma]$ is in
the kernel of~$h$,
so $[\Gamma,\Gamma]$ sends every leaf to itself. Furthermore,
Lemma~\ref{finitecommutator} below asserts that
$[\Gamma,\Gamma]$ is a finite-index subgroup of~$\Gamma$.
 \end{pf}

\begin{lem} \label{finitecommutator}
 If each homomorphism of a group  $\Gamma$ into the group of homeomorphisms
of the circle has a finite
orbit,  then   $A:=\Gamma/[\Gamma,\Gamma]$ is a finite group.
 \end{lem}

\begin{pf}
 $A$ is a finitely generated abelian group and if it is not finite we can
find a homomorphism from
$\Gamma$ onto $\Bbb Z$; but $\Bbb Z$ clearly  acts faithfully on the circle
with no
finite orbits (for example, by mapping $1\in \Bbb Z$ to
an irrational rotation), so the action of  $\Gamma$ would not have any
finite orbits.
 \end{pf}

\section{Proof of Theorem~\ref{bigtwo}}

The following lemma is a slight generalization of the Thurston Stability
Theorem \cite{thurston}:
 \begin{itemize}
 \item we only assume the action of~$\Gamma$ is germinal, rather than being
globally defined; and
 \item we do not assume that each element of~$\Gamma$ is defined in a
neighborhood of~$0$, but only
on a set~$X$ that accumulates at~$0$.
 \end{itemize}
 For completeness, we provide a proof, although Thurston's original proof
can easy be generalized to
this setting.

\begin{lem}[{{\rm Thurston, cf.\ \cite{thurston}}}] \label{thurston}
Suppose
 $\Gamma$~is a finitely generated group,
 $X$~is a compact subset of $[0,1]$ that accumulates at~$0$,
 and, for each $\gamma \in \Gamma$, we have a
$C^1$ diffeomorphism $\phi_\gamma \colon [0,a_\gamma)\rightarrow
[0,b_\gamma)$, for positive constants
$a_\gamma$ and~$b_\gamma$.
 Assume
\begin{enumerate}
\item \label{Thurston-comm} $\Gamma/[\Gamma,\Gamma]$ is finite; and
\item for each $\gamma_1,\gamma_2\in \Gamma$, there exists $c>0$, such that
 $\phi_{\gamma_1 \gamma_2}|_{[0,c]\cap
X}=\phi_{\gamma_1}\phi_{\gamma_2}|_{[0,c]\cap X}$.
\end{enumerate}
 Then there exists $a > 0$, such that, for every $\gamma \in \Gamma$ and $x
\in [0,a) \cap X$, we have
$\phi_\gamma(x) = x$.
\end{lem}

\begin{pf} (cf.\ \cite{reeb}) We use nonstandard analysis. Let $x\in
\widehat{X}$
with $x\thickapprox 0$. (We use
$\widehat{X}$ to denote the nonstandard set corresponding to $X$, and we write
$a\thickapprox b$ if $a-b$ is infinitesimal.) It suffices to show that
$\phi_\gamma(x) = x$ for every $\gamma \in \Gamma$. Suppose not, and let
$\epsilon = \text{max}_{\gamma\in F}|\phi_\gamma(x) - x|$,
where $F$ is some (fixed) finite generating set for~$\Gamma$.

Define $d \colon \Gamma \to (0,\infty)$ by
 $$d(\gamma) = \text{Re} \left( \frac{\phi_\gamma(x)}{x} \right) =
\phi_\gamma'(0).$$
 Then $d$ is a homomorphism, so (\ref{Thurston-comm}) implies $d(\gamma) =
1$ for every~$\gamma$.
 Therefore, for each $\gamma \in \Gamma$, we have
 \begin{equation} \label{Thurston-x+o(x)}
 \mbox{$\phi_\gamma(y) = y + o(y)$ whenever $y \approx 0$.}
 \end{equation}

For $\gamma \in \Gamma$, and $y_1,y_2 \in \widehat X$ with $y_1 \approx y_2
\approx 0$,
the Mean Value Theorem  implies there is some $\xi \approx 0$, such that
 $$ \phi_\gamma(y_1) - \phi_\gamma(y_2)
 = \phi_\gamma'(\xi) (y_1 - y_2) .$$
 Furthermore, because $\phi_\gamma'(0) = 1$ and $\phi_\gamma'$ is
continuous, we have
$\phi_\gamma'(\xi) = 1 + o(1)$, so
  \begin{equation} \label{Thurston-x-y}
 \phi_\gamma(y_1) - \phi_\gamma(y_2)
 = \bigl(1 + o(1) \bigr) (y_1 - y_2)
 = y_1 - y_2 + o(y_1 - y_2)
 . \end{equation}

By induction on word
length (and using Eq.~(\ref{Thurston-x-y}), much as in the multi-line
displayed equation below), we see that
$\phi_\gamma(x) - x =O(\epsilon)$, for every $\gamma\in \Gamma$. Thus, we
may define
$f\colon \Gamma\rightarrow \Bbb R$ by
 $$ f(\gamma)=\text{Re}\left(\frac{\phi_\gamma(x) - x}{\epsilon}\right).$$
 For $\gamma,\lambda\in  \Gamma$, we have
\begin{align*}
f(\gamma \lambda)&\thickapprox \frac{\phi_{\gamma\lambda}(x) - x}{\epsilon}\\
&=\frac{\phi_\gamma(x) - x}{\epsilon} +\frac{\phi_{\gamma\lambda}(x) -
\phi_\gamma(x)}{\epsilon}\\
&= \frac{\phi_\gamma(x) - x}{\epsilon} + \frac{\phi_\lambda(x)-x +o \bigl(
\phi_\lambda(x)-x \bigr)}{\epsilon}\\
&\thickapprox \frac{\phi_\gamma(x) - x}{\epsilon}
+\frac{\phi_\lambda(x)-x}{\epsilon}\\
&\thickapprox f(\gamma) +f(\lambda),
\end{align*}
so $f(\gamma\lambda)=f(\gamma) + f(\lambda)$, which means that $f$ is a
homomorphism.
Because $\Gamma/[\Gamma,\Gamma]$ is finite, we conclude that $f(\gamma)=0$ for
every $\gamma \in \Gamma$. Therefore $\phi_\gamma(x) = x$ for every
$\gamma\in \Gamma$.
 \end{pf}

All the assumptions of Theorem \ref{bigtwo} are in force from now on.

\begin{lem}\label{finalone}
 Let $\Gamma'$ be a finite-index subgroup of~$\Gamma$ and suppose that $W$ is a
connected, open, $\Gamma'$-invariant, ${\cal F}$-saturated subset of~$M$
whose boundary components
are $\Gamma'$-invariant compact leaves. If the compact leaves accumulate on
a boundary component~$L$
of $\widehat{W}$, then all the compact leaves in some foliated neighborhood
of $L$ in $\widehat{W}$
are also $\Gamma'$-invariant.
 \end{lem}

\begin{pf}
 Let $U$ be a connected, $\widehat{\cal F}$-saturated neighborhood of~$L$
in~$\widehat{W}$. We may assume $U$ is so
small that there is a complete $C^1$ transversal $\alpha\colon [0,1] \to U$
of~$U$, such that each compact leaf in~$U$
meets~$\alpha$ exactly once. We may assume $\alpha(0) \in L$.

As the $\Gamma$-action and the foliation (hence the local holonomy maps)
are $C^1$,  we can construct, for each
$\gamma\in \Gamma'$,  a $C^1$ diffeomorphism $\phi_\gamma \colon
[0,a_1]\rightarrow [0,a_2]$, for positive constants
$a_i$, such that $\alpha \bigl( \phi_\gamma(t) \bigr)$ is on the same leaf
as $\gamma\bigl( \alpha(t) \bigr)$,
for each $t \in [0,a_1]$.
 Let
 $$X = \{\, x \in [0,1] \mid \mbox{$\alpha(x)$ is on a compact leaf} \,\} .$$
 From Lemma~\ref{finitecommutator}, we know that
$\Gamma'/[\Gamma',\Gamma']$ is finite.  Therefore,
Lemma~\ref{thurston} applies, so we conclude that there is an interval
$[0,a)$, such that $\alpha(x)$ is on the
same leaf as~$\gamma\bigl( \alpha(x) \bigr)$, for all $x \in X \cap [0,a)$.
So $\Gamma'$ fixes the compact leaf
containing $\alpha(x)$, for each $x \in X \cap [0,a)$.
 \end{pf}

\begin{lem}\label{aux}
If a finite group $H$ acts on the circle by orientation preserving
$C^1$-diffeomorphisms, then the
action of the commutator subgroup $[H,H]$ is trivial.
\end{lem}

\begin{pf}
By a standard averaging procedure $H$ can be assumed to act by isometries
of a Riemannian metric on the circle.
Therefore the action corresponds to a homomorphism of $H$ into the group of
rotations of the circle, which
is an abelian group. If follows that the commutator subgroup must act
trivially.
\end{pf}

\begin{lem}\label{nocpct->done}
 Let $\Gamma'$ be a finite-index subgroup of~$\Gamma$ and suppose that
  $W\subset M$ is a connected, open, $\Gamma'$-invariant, $\cal
F$-saturated subset, such that the boundary
components of $W$ are $\Gamma'$-invariant compact leaves, and every compact
leaf in~$W$ is $\Gamma'$-invariant. Then
all leaves in $W$ are $[\Gamma',\Gamma']$-invariant.
\end{lem}

 \begin{pf}
 The complement in $W$ of the set of
compact leaves is a countable union of connected, $\Gamma'$-invariant,
$\cal F$-saturated sets $U_i$, $i=1,2, \dots$,
such that the boundary of each $U_i$ is the union of at most $2$ compact,
$\Gamma'$-invariant leaves.
Since there are no compact leaves in each of the $U_i$, each leaf of ${\cal
F}|_{U_i}$ is, by assumption, without
holonomy. Theorem \ref{general-two} implies that either
 \begin{enumerate}
 \item \label{fiber-circle} ${\cal F}|_{U_i}$ fibers over the circle; or
 \item \label{fiber-I}  $(\widehat U_i, \widehat{\cal F})$ fibers over
$[0,1]$; or
 \item \label{fiber-denseleaf} every leaf
of  ${\cal F}|_{U_i}$ is dense in $U_i$.
 \end{enumerate}

In Case~\ref{fiber-circle}, we obtain a $C^1$ action of $\Gamma'$ on the
circle.
This action has to be finite, that is, a finite-index subgroup of~$\Gamma'$
must
leave invariant each leaf in $U_i$. From Lemma~\ref{aux}, we know that this
normal subgroup contains the commutator
subgroup $[\Gamma',\Gamma']$.

In Case~\ref{fiber-I}, we obtain a $C^1$ action of $\Gamma'$ on $[0,1]$.
Lemma~\ref{thurston} implies that this action is trivial.

In Case~\ref{fiber-denseleaf}, we can apply
Theorem~\ref{general-three}(\ref{partone}) and
Lemma~\ref{corobigtwo}(\ref{corobigtwo-leaf}) to conclude that each leaf in
$U_i$ is
$\Gamma'$-invariant.

Therefore, $[\Gamma',\Gamma']$ fixes every leaf of $W$.
\end{pf}

{\it We now conclude the proof of Theorem \ref{bigtwo}.} It suffices to
show that there is a
finite-index subgroup~$\Gamma'$ of~$\Gamma$, such that $\Gamma'$ fixes
every compact
leaf, for then Lemma~\ref{nocpct->done} shows that $[\Gamma',\Gamma']$
fixes every leaf (and
Lemma~\ref{finitecommutator} implies that $[\Gamma',\Gamma']$ has finite
index in~$\Gamma$).

We may assume that ${\cal F}$ has a compact leaf. Then
Theorem~\ref{atomthm} states that
there is a finite-index subgroup~$\Gamma'$ of~$\Gamma$, and a compact
leaf~$L_0$
of~${\cal F}$, such that $\Gamma'$ fixes~$L_0$. Let $W_0$ be the complement
of the union of all the
compact leaves of~${\cal F}$.  From Theorem~\ref{lemgen}, we
know that only finitely many components $U_1,U_2,\ldots,U_k$ of~$W_0$ fail
to be foliated products,
and that each of these components has only finitely many boundary leaves, so,
replacing $\Gamma'$ with a finite-index subgroup, we may assume that
 \begin{equation} \label{fixthebads}
 \mbox{$\Gamma'$ fixes each $U_i$, and its boundary components.}
 \end{equation}

We claim that $\Gamma'$ fixes every compact leaf. (This will complete the
proof.)

Because $M$ is connected, it suffices to show, for each
$\Gamma'$-invariant compact leaf~$L$, that there are connected, open,
$\Gamma'$-invariant, ${\cal F}$-saturated subsets $U^+$ and~$U^-$ of~$M$
whose boundary consists of
compact leaves, such that
 \begin{itemize}
 \item $\Gamma'$ fixes each compact leaf in the closure of~$U^+ \cup U^-$;
 \item $L$ is a border leaf of both~$U^+$ and~$U^-$; and
 \item $U^+$ is on the positive side of~$L$ and $U^-$ is on the negative
side of~$L$.
 \end{itemize}

Fix a $\Gamma'$-invariant compact leaf~$L$.

If the compact leaves accumulate on the positive side of~$L$, then
Lemma~\ref{finalone} provides an appropriate $\Gamma'$-invariant open
set~$U^+$.

If the compact leaves do not accumulate on the positive side of~$L$, then,
on its positive side,
$L$ is part of the boundary of
a component~$U^+$ of~$W_0$. If $U^+$ is not a foliated product, then, from
\eqref{fixthebads},
we know that $\Gamma'$ fixes each boundary component of~$U^+$, as desired.
Now suppose that $U$ is a foliated product.
Then the boundary of~$U$ consists of at most two compact leaves. One of
these boundary leaves is $L$,
which is known to be $\Gamma'$-invariant. Thus the second boundary leaf (if
it exists) must also be $\Gamma'$-invariant.

The construction of~$U^-$ is similar. \qed

\end{document}